\documentclass[12pt]{article}
\usepackage{amssymb,latexsym,amsmath}
\def\C{\mathbb C}
\def\R{\mathbb R}
\def\N{\mathbb N}
\def\Z{\mathbb Z}

\DeclareMathOperator{\cosec}{cosec}

\newtheorem{thm}{Theorem}[section]
\newtheorem{prop}{Proposition}[section]
\newtheorem{lem}{Lemma}[section]
\newtheorem{con}{Conjecture}[section]
\newtheorem{definition}{Definition}[section]

\begin{document}
\sffamily
\title{Non-real zeros of  derivatives}
\author{J.K. Langley}
\maketitle

\centerline{\textit{Dedicated to the memory of Larry Zalcman} }

\begin{abstract}
A number of results are proved concerning non-real zeros of 
derivatives of real meromorphic
functions. 
\end{abstract}

\section{Introduction}\label{intro}

This paper concerns non-real zeros of derivatives
of real meromorphic functions in the plane, that is, meromorphic functions mapping
$\R$ into $\R \cup \{ \infty \}$. 
The case  of real entire functions has seen extensive research
\cite{ BEpolya, BEL,CCS2, HelW1, KiKim,  lajda, LeO, SS},
motivated at least in part by the Wiman conjecture  
(proved in \cite{BEL,LeO,SS})
that if $f$ is a real entire function and $f$
and $f''$ have only real zeros, then $f$ belongs to
the Laguerre-P\'olya class consisting of
locally uniform limits of real polynomials
with real zeros.
The following theorem combines results from \cite{BEL,Larossi}.

\begin{thm}[\cite{BEL,Larossi}]
 \label{thm1}
 Let $f$ be a real meromorphic function of infinite order in the plane 
such that $f$ or $1/f$ has finitely many poles and non-real zeros.
Then $f''/f'$ has infinitely many non-real zeros, that is, $f''$ has infinitely many non-real
zeros which are not zeros of $f'$. 
\end{thm}

This result
will be strengthened as follows. 

\begin{thm}
 \label{thm1A}
 Let $f$ be a real meromorphic function of infinite order in the plane,
 with finitely many non-real zeros and poles. 
 Assume that $f = f_1/f_2$, where the $f_j$ are real entire functions, with no common zeros,
 and that at least one of $f_1$ and $f_2$ has finite lower order. 
Then $f''/f'$ has infinitely many non-real zeros.
\end{thm}

Thoerem \ref{thm1A} represents a fairly substantial improvement of 
Theorem \ref{thm1}, since under the hypotheses of the latter one of $f_1, f_2$ may be
assumed to be a polynomial. The theorem is applicable, in particular,  if $f$ is a
real meromorphic function of infinite order in the plane,
 with finitely many non-real zeros and poles, for which the exponent of convergence 
 of either the zeros or the poles of $f$ is finite.
Theorem \ref{thm1A} will be deduced from the next result: here, and subsequently,
$H$ denotes the open upper half-plane $\{ z \in \C : \, {\rm Im} \, z > 0 \}$. 

\begin{thm}
\label{thm1AA}
Let $L$ be a real meromorphic function in the plane, with finitely many non-real poles, 
and assume that $L$ has a representation 
\begin{equation}
 L = h R_1 \psi_1 + R_2 \psi_2 ,
\label{Lnewrep}
\end{equation}
in which:
$h$ is a real transcendental entire function;
the $R_j$ are real rational functions, with $R_1 \not \equiv 0$;
each $\psi_j$ satisfies either $\psi_j \equiv 1$ or $\psi_j(H) \subseteq H$.

Then $L + L'/L$ has infinitely many non-real zeros. 
\end{thm}

To deduce Theorem~\ref{thm1A} from Theorem~\ref{thm1AA} it will suffice to show that
$L = f'/f$ has a representation (\ref{Lnewrep}):
this follows from the formula $L + L'/L = f''/f'$. Such a representation is well known if  $f$ is as in Theorem \ref{thm1}, because   $f'/f$ then has a
 \textit{Levin-Ostrovskii factorisation} of form (\ref{Lnewrep}) with $R_2 \equiv 0$
 \cite{BEL,Larossi,LeO} (see also Lemma \ref{lemlevost}). 


The methods of Theorems \ref{thm1A} and \ref{thm1AA} turn out also to be
applicable to 
a strand initiated in \cite{Lajda09}, in which the second derivative
$f''$ is replaced 
by $f'' + \omega f$, with $\omega \in \R \setminus \{ 0 \}$. 
Here attention is necessarily restricted to the case $\omega > 0$, 
as illustrated by an example 
cited in \cite{Lajda09}: for $a \in \R$ writing  $f'(z)/f(z) = a + e^{-2az} $  makes
$f$, $1/f$ and  $f''(z) - a^2 f(z) =  e^{-4az} f(z)$ all zero-free.

\begin{thm}
 \label{cor1}
 Let $f$ be as in the hypotheses of Theorem \ref{thm1A},
 and let $\omega $ be a positive real number. Then $f'' + \omega f$ has infinitely many non-real zeros. 
\end{thm}

When $f$ is a real entire function of infinite order with finitely many non-real zeros,
Theorem~\ref{cor1} is not new \cite{Lajda09}, but the present proof 
is considerably simpler than that of \cite{Lajda09} and the result substantially more general.
For  results on non-real zeros of $f''+\omega f$ 
when $\omega \geq 0$ and $f$ has finite order, the reader is referred to \cite{Lajda09,Lawiman10,Lawiman13} and  \cite[Theorem 1.5]{Laldp13}. 

As with the proof of Theorem \ref{thm1A}, Theorem \ref{cor1} will be deduced from a result
involving functions of the form (\ref{Lnewrep}).

\begin{thm}
 \label{thmlindiffpol}
Let $L$ be as in the hypotheses of Theorem \ref{thm1AA},
and let $a, b$ be  positive real numbers. 
Then $L' + a L^2+ b$ has infinitely many non-real zeros.
\end{thm}
The simple example $L(z)= \tan z$, $a = b = 1$, shows that the 
requirement that $h$ is transcendental in (\ref{Lnewrep}) is not redundant in Theorem \ref{thmlindiffpol}.
The following property will play a pivotal role in the
proofs of Theorems \ref{thm1AA} and \ref{thmlindiffpol}.

\begin{definition}
\label{defnUHWV}
A transcendental meromorphic function $L$ in the plane has the UHWV 
property if  there exist
$\tau, \gamma$ with
$1/2 < \tau < \gamma < 1$, an unbounded subset $E_1 $ of $[1, + \infty )$ 
and a function $N(r): E_1 \to (1, + \infty )$ satisfying the following:\\
(A) $\lim_{r \to + \infty, r \in E_1} N(r) = + \infty $;\\
(B) for each 
$r \in  E_1$ there exists $z_0 = z_0(r)$ with 
\begin{equation}
\label{z0def}
|z_0| = r, \quad N(r)^{-\tau} <  \arg z_0 < \pi -  N(r)^{-\tau} ,
\end{equation}
such that, uniformly as $r \to + \infty $ in $E_1$,
\begin{equation}
\label{wv1a}
L (z) \sim L(z_0) \left( \frac{z}{z_0} \right)^{N(r)} \quad 
\hbox{and} \quad N(r)^{1/2} = o \left(   \log^+ |L(z)| \right) 
\end{equation}
on 
\begin{equation}
\label{Qrdef}
Q_r = \left\{ z \in \C : \, 
 \left| \log \frac{z}{z_0} \right| \leq N(r)^{-\gamma} \right\} .
\end{equation}
\end{definition}

Here UHWV stands for ''upper half-plane Wiman-Valiron'' and
standard results from the Wiman-Valiron theory \cite{Hay5} imply
 that if $L$ is a real transcendental entire function then $L$ has the UHWV property, 
 with $N(r)$ the central index (see Lemma \ref{auxwvlem}).

\begin{thm} 
\label{thmUHWV1}
Let $L$ be a real transcendental meromorphic function in the plane with 
finitely many non-real poles and assume that $L$ has the UHWV 
property. Then $L+L'/L$ has infinitely many non-real zeros. 
\end{thm}

To prove Theorem \ref{thm1AA} it will suffice to show that $L$ has the UHWV property and 
apply
Theorem \ref{thmUHWV1} 
 directly. It is not clear whether, under the hypotheses of Theorem \ref{thmUHWV1}, 
 $L'+aL^2 + b$ automatically has non-real zeros when $a, b > 0$, and the proof of
 Theorem \ref{thmlindiffpol} will use the representation (\ref{Lnewrep}) alongside
 the UHWV property.

The next focus of the present paper is
 the general problem of classifying all real
meromorphic functions in the plane which, together with some of their derivatives,
have only real zeros and poles 
\cite{HSW,Hin1,Hin3}. In this direction, the following conjecture was advanced in 
\cite{HSW}. 

\begin{con}[\cite{HSW}]
 \label{conhsw}
Let $f$ be a real transcendental meromorphic function in the plane with at least one pole, and assume that all zeros and poles of
$f$, $f'$ and $f''$ are real, and that all poles of $f$ are simple. Then 
$f$ satisfies
\begin{equation}
 \label{hswform}
f(z) =  C  \tan (az+b) + Dz + E, \quad  a, b, C, D, E \in \R. 
\end{equation}
\end{con}
If $f$ is allowed multiple poles then there are further examples 
for which $f$, $f'$ and $f''$ have only real zeros and poles
 \cite{Hin4}. Results  from 
\cite{HSW,HinRos,Laams09,ankara,Lawiman13,nicks} show  that the conjecture is true  if, in 
addition, $f'$ omits some finite value. 
Furthermore,  Theorem \ref{thm1A} and \cite{Lawiman10}
together show that there are no
functions $f$ satisfying the hypotheses of 
Conjecture \ref{conhsw} such that either of the following holds: 
$f$ has infinite order and the zeros or poles of $f$ have finite exponent of convergence;
$f$ has finite order and infinitely many poles but finitely many zeros. 
The conjecture 
was also proved in \cite[Theorem 1.4]{Lajda18} for 
 real transcendental meromorphic functions in the plane
 which map the open upper half-plane $H$ 
into itself. All zeros and poles of such functions are automatically  real and simple
and \textit{interlaced}   \cite{Le}: that is,  
between any two consecutive poles of $f$ there is a zero, and 
between consecutive zeros of $f$ lies a pole (this follows from a consideration of residues for
$f$ and $1/f$).


\begin{thm}
\label{cor2}
Let $f$ be a real transcendental meromorphic function in the plane with 
infinitely many zeros and poles, all real, simple and interlaced. 

If $f''$ has finitely many non-real zeros then
\begin{equation}
\label{jdaform}
f(z) = D \left( Az+ B + \frac{R(z)e^{icz}-1}{A_1R(z)e^{icz} - \overline{A_1} } \right)  ,
\end{equation}
where $A, B, c, D$ and $A_1$ are constants with $A, B, c, D$ real and $A_1 \in H$, 
while $R$ is a rational function with  $|R(x)| = 1$ for all real $x$.

If $f''$ has only real zeros then $f$ is as in (\ref{hswform}). 
In particular, 
 Conjecture \ref{conhsw} is true under the additional hypothesis that $f$ has infinitely many zeros and poles, all simple and interlaced. 
 \end{thm}

Theorem \ref{cor2} will be deduced from
\cite[Theorem 1.4]{Lajda18}  and the
 following result involving real meromorphic functions with real zeros and poles such that,
with finitely many exceptions, all poles are simple and 
adjacent poles are separated by at least one zero.

 \begin{thm}
 \label{thm11}
 Let $U$ be a real meromorphic function 
 in the plane,
 all but finitely many of whose zeros are real.
 Assume further that $U$ has infinitely many poles $X$, 
 all but finitely many of which  are 
 real and simple 
 and have a corresponding real
 zero $Y $ of $U$ with $X < Y$ and 
 $U(x) \neq \infty$ on $(X, Y)$. 
 Then $U$ satisfies the following.\\
 (i) $U$ has a representation 
 \begin{equation}
 U = S \psi ,
 \label{ww1}
 \end{equation} 
 where $S$ is a real meromorphic function in the plane with finitely many poles and 
 $\psi(H) \subseteq H$.\\
 (ii) If  $S$ has infinite order then $U''/U'$ has infinitely many non-real zeros.\\
 (iii) If  $S = Re^P$, 
with $R$ a real rational function and  $P$  a non-constant polynomial,
 then $U^{(m)} $ has infinitely many non-real zeros, for each $m \geq 2$. 
 \end{thm}
Part (i) is not new, but its inclusion  is convenient for the statement
 and proof of 
 parts (ii) and (iii). The standard construction
 is outlined in Lemma~\ref{lemlevost}: here  $\psi$ is determined only up to a rational factor, 
 but the choice of such a factor does not affect (ii) or (iii).
 
  If $S$ is transcendental
with finitely many zeros and poles in (\ref{ww1}), 
then either (ii) or (iii) is applicable, 
although 
Theorem \ref{thm11} says nothing about
 the case where $S$ has finite order and infinitely many zeros.
Simple examples such as $\cot z$ 
show that Theorem \ref{thm11}(iii) fails  for $m=2$ if $P$ is constant: 
for an example not of the form (\ref{hswform}) set
$$
V(z) = z \cot z, \quad 
V''(z) = -2 \cosec^2 z + 2z \cosec^2 z \cot z 
= 2 ( z - \tan z ) \cosec^2 z \cot z . 
$$
Since the iterates of $\tan z$ converge to $0$ on $\C \setminus \R$, 
all fixpoints of $\tan z$ are real, and so are all zeros of $V''$.

This paper is organised as follows. After preliminary considerations in Sections 
\ref{prelim} and \ref{singu}, Theorem \ref{thmUHWV1} is proved in 
Section \ref{pfuhwv1}. The UHWV property is discussed further in Section \ref{sufficient},
the proofs of 
Theorems \ref{thm1A} and \ref{thm1AA} 
appearing in Section \ref{pfthm1AA}. 
Next, Theorem \ref{thmlindiffpol} is proved, and Theorem \ref{cor1} is deduced from it,
in Section \ref{pfthmlindiffpol}. Finally,
Theorem \ref{thm11} is established in  Section  \ref{infiniteorder},
and Theorem \ref{cor2}  
in Section \ref{pfcor2}.

\section{Preliminary lemmas} \label{prelim}

\begin{lem}
\label{lemlevin}
There exists a positive constant $c_0$
such that if $\psi : H \to H$ is analytic then, for 
 $ r \geq 1$ and 
$\theta\in (0,\pi)$,
\begin{equation}
\frac{|\psi(i)| \sin\theta}{5 r} < |\psi(re^{i\theta})| <
\frac{5r |\psi(i)|}{\sin\theta}\quad \hbox{and}  \quad 
 \left| \frac{ \psi'( re^{i \theta } ) }{ \psi( re^{i \theta } ) } \right|
 \leq \frac{c_0}{r \sin \theta} .
  \quad\label{C1}
\end{equation}
\end{lem}
Both of these estimates  are standard: the first is essentially just Schwarz' lemma
 \cite[Ch. I.6, Thm $8'$]{Le}, while the second follows from Bloch's theorem applied to $\log \psi$.
\hfill$\Box$
\vspace{0.1in}

\begin{lem}
[\cite{EC1}]
 \label{thmel}
Let  $\Omega$ be  a plane domain. 
Let $\mathcal{L}$ be the family of all analytic functions  $L$ on $\Omega$ such that 
$\Psi_2(L) + 1 = L' + L^2 + 1$ has no zeros on $\Omega$. 
Then $\mathcal{L}$ is normal. 
\end{lem}
Lemma \ref{thmel} is a special case of  \cite[Theorem 4]{EC1}, which was proved by  Eleanor Lingham (under her maiden name Clifford), using the highly effective and influential
rescaling technique invented by Larry Zalcman and developed further by him and Pang Xuecheng \cite{Zalc}. 
\hfill$\Box$
\vspace{0.1in}

\begin{lem}[\cite{EF2}]\label{lemEF}
Let $1 < r < R < + \infty $ and let the function $g$ be meromorphic in
$|z| \leq R$. Let $I(r)$ be a subset of $[0, 2 \pi ]$ of Lebesgue
measure $\mu (r)$. Then
$$
\frac1{2 \pi} \int_{I(r)} 
\log^+ | g(re^{i \theta })| \, d \theta \leq 
\frac{ 11 R \mu (r)}{R - r} \left( 1 + \log^+ \frac1{\mu(r)} \right)
T(R, g).
$$
\end{lem}

\begin{lem}[\cite{Hay3}]\label{Haylem}
Let 
$S(r)$ be an unbounded positive
function on $[1, + \infty )$ which is  non-decreasing and continuous from the right.
Let $A > 1, B > 1$ and $G = \{ r \geq 1 : S(Ar) \geq B S(r) \}$. Then
the upper logarithmic density of $G$ satisfies
$$
\overline{{\rm logdens }} \, G = \limsup_{r \to \infty} \left( \frac1{\log r} \, \int_{[1, r] \cap G} \, \frac1t \, dt \right) 
\leq  \left( \frac{\log A}{\log B} \right) \limsup_{r \to \infty} \frac{ \log^+ S(r)}{\log r } . 
$$
\end{lem}

The next lemma proves Theorem \ref{thm11}(i). 

 \begin{lem}
 \label{lemlevost}
 Let $U$ be as in the hypotheses of Theorem \ref{thm11}(i). 
Then $U$ has a representation (\ref{ww1}) 
in which $S$ is a real meromorphic function with finitely many poles
in the plane and  $\psi(H) \subseteq H$. 
\end{lem}
\textit{Proof.}
This is the standard Levin-Ostrovskii construction \cite{BEL,LeO}. 
By assumption, all but finitely many poles of $U$ are real and simple, and all but finitely many of
these may be 
labelled 
$x_k$
in such a way  that $x_k < x_{k+1}$
and there is a zero $y_k$ of $U$ with $x_k < y_k <  x_{k+1}$ and $y_k/x_k > 0$. 
If $K$ denotes the set of these $k$, the product 
\begin{equation}
 \label{wwpsi}
 \psi(z) = \prod_{k \in K} \frac{1-z/y_k}{1-z/x_k} 
 \end{equation}
 converges by the alternating series test, and  maps $H$ into $ H$
 because, for $z \in H$,
 $$
 \arg \psi (z) 
  =  \sum_{k \in K}  \arg \frac{y_k-z}{x_k-z} \in (0, \pi ). 
 $$
 This proves Lemma \ref{lemlevost}, with $\psi$ determined up to a rational factor. 
  \hfill$\Box$
\vspace{0.1in} 
 
 \begin{lem}
\label{lem2}
Let  $m \geq 2$ be an even integer
and let 
$$
Q_m(y) =  \sum_{j=0}^{m}    { m \choose j } (j+1)! \, y^{j} .
$$
Then there exists $d_m > 0$ such that 
$Q_m(y) \geq d_m \max \{ 1, \, y^{m} \} $
for all $y \in \R$. 
\end{lem}
\textit{Proof.} 
Since 
$Q_m(0) = 1$ it suffices to show that $Q_m(y) > 0$ for $y \neq 0$,
this
being obvious if $y > 0$. 
 For $y < 0$  write $x = 1/y < 0$ and $ P(x) = e^{-x} x^{-2} $.  
 Then  Leibniz' rule and the fact that $m$ is even together yield
\begin{eqnarray*}
 P(x) &=& \frac1{x^2} - \frac1x + \frac1{2!} - \frac{x}{3!} + \ldots + \frac{x^m}{(m+2)!} -\frac{x^{m+1}}{(m+3)!}  + \dots,\\
 0 &<& 
 \frac{(m+1)!}{x^{m+2}} - \frac{m!} {x^{m+1}}
 +   \frac{m!}{(m+2)!}  - \frac{(m+1)!}{(m+3)!} \, x + \dots \\
 &=&  P^{(m)}(x) 
= \sum_{j=0}^{m}    { m \choose j } (-1)^{m-j} e^{-x} (-1)^j (j+1)! \, x^{-2-j}\\
&=& \sum_{j=0}^{m}    { m \choose j } e^{-x} (j+1)! \, x^{-2-j} =  x^{-2} e^{-x}  Q_m(y) .
\end{eqnarray*}
\hfill$\Box$
\vspace{.1in}

\section{Transcendental singularities of the inverse function}\label{singu}

Throughout this section let $G$ be a transcendental meromorphic 
function in the plane.
Suppose first that  $G(z) \to a \in \C \cup \{ \infty \}$ as $z \to \infty$ along a path
$\gamma $; then the inverse $G^{-1}$ is
said to have a transcendental singularity
over the asymptotic value $a$~\cite{BE,Nev}. If $a \in \C$ then for each $\varepsilon > 0$ 
there  exists a component $\Omega = \Omega( a, \varepsilon, G)$ of the
set $\{ z \in \C : |G(z) - a | < \varepsilon \}$ such that
$\gamma \setminus \Omega$ is bounded, these components being called neighbourhoods of the singularity \cite{BE}.
Two  paths $\gamma, \gamma'$ on which $G(z) \to a$ determine distinct singularities if the corresponding components
$\Omega( a, \varepsilon, G)$, $\Omega'( a, \varepsilon, G)$ are disjoint for  some $ \varepsilon > 0$.
The singularity 
is called direct  \cite{BE} if $\Omega( a, \varepsilon, G)$, for some
$\varepsilon > 0$, contains finitely many zeros of $G-a$, and 
indirect otherwise.
A transcendental singularity will be referred to as lying in an open set $D$ 
if $\Omega(a, \varepsilon, G) \subseteq D$ for all sufficiently small positive
$\varepsilon $. Transcendental singularities over $\infty$ may be 
classified using $1/G$.

The following lemmas from \cite{Lajda09,Lawiman13} 
link asymptotic values approached on paths in $H$ with the growth of 
the Tsuji characteristic $ \mathfrak{T} (r, g) $ \cite{BEL,GO,Tsuji}, which is 
defined for meromorphic functions $g$ 
on the closed upper half-plane $\overline{H} =
\{ z \in \C : \mathrm{Im} \, z \geq  0 \}$.

\begin{lem}[\cite{Lawiman13}, Lemma 2.2]
\label{directsinglem2}
Let $L \not \equiv 0$ be a real  meromorphic function in the plane such that
$\mathfrak{T}(r, L) = O( \log r )$ as $r \to \infty$, and define $F$ by  $F(z) = z-1/L(z)$. 
Assume that at least one of $L$ and $1/L$ has finitely many non-real poles. 
Then there exist finitely many $\alpha \in \C $ such that $F(z)$ or $L(z)$
tends to $\alpha$ as $z$ tends to infinity along a path in $\C \setminus \R$.
\end{lem}

\begin{lem}[\cite{Lajda09}, Lemma 2.4]\label{directsinglem}
Let $G$ be a meromorphic function in the plane such that
$\mathfrak{T}(r, G) = O( \log r )$ as $r \to \infty$. Then there is 
at most one direct transcendental singularity of $G^{-1}$ lying in 
$H$.
\end{lem}

The following proposition is a stronger version of \cite[Lemma 3.2]{Lajda09}, with a simpler proof,
which will occupy the remainder of this section. 
Here
$B(a, r)$ denotes the open ball of centre $a \in \C$ and radius $r > 0$.

\begin{prop}
 \label{propnumberzeros}
Suppose that $R \in (0, + \infty)$ and the transcendental meromorphic function
$G$ has no asymptotic values $w$ with $0 < |w| < R  < \infty  $, and finitely many critical values $w$ with $|w| < R$. 
Let $A$ be a component of the set $G^{-1}(B(0, R))$. Then the number of  zeros of $G$ in $A$, counting multiplicities, plus the number 
of transcendental singularities of $G^{-1}$ over $0$, lying in $A$, exceeds by at most
$1$ 
the  number  of zeros of $G'$ in $A$, again counting multiplicities.
\end{prop}
\textit{Proof.} 
It may be assumed that 
there  exists a component $A$ of  $G^{-1}(B(0, R))$ which contains a finite number, $M$ say,
of zeros of $G'$,
counting multiplicities, but also contains zeros $u_1, \ldots, u_p$ of $G$, repeated according to multiplicity, as well as 
$q$  pairwise disjoint neighbourhoods $\Omega_j(0, s, G)$ of transcendental singularities of $G^{-1}$ over $0$, where $s > 0$ is small and $M+1 \leq p + q < \infty$.
It is not assumed at this stage that there are no other zeros of $G$, nor other  transcendental singularities of $G^{-1}$ over $0$, lying in $A$, nor even that the number
of these is finite. 
Choose points $v_j \in \Omega_j(0, s, G)$, for $j=1, \ldots, q$. Then 
$u_1, \ldots, u_p, v_1 , \ldots, v_q$ may be joined to each other by paths in $A$
and so all lie in a compact connected subset of $A$ on which $|G(z)| \leq S_1$, and hence in a
component $B \subseteq A$ of  $G^{-1}(B(0,S_2))$, for some $S_1, S_2$ with $s < S_1 < S_2 < R$. 

These observations show that  it is enough to prove that $p+q \leq M+1$ when
$G$ has no critical or asymptotic values $w$ with $|w| = R$. 
Let $w_1, \ldots, w_N$ be the critical values of $G$ with $0 < |w| < R$. Join each $w_j$ to a point $w_j^*$ on $|w| = R$ by a 
straight line segment $\lambda_j$ 
in the annulus $2s  < |w| \leq R$, in such a way that these $\lambda_j$ are pairwise disjoint; if the $w_j$ have distinct arguments modulo $2 \pi$, the 
$\lambda_j$ may be taken to be radial segments, while if repetition occurs the segments may be rotated slightly about $w_j$. 
Let $E_0 = B(0, R)$ and, for $m=1, \ldots, N$, set
$$
E_m = E_{m-1} \setminus \lambda_m = 
E_0 \setminus \left( \bigcup_{j=1}^m \lambda_j \right) .
$$
Since $E_N \setminus \{ 0 \}$ contains no asymptotic or critical values of $G$, a straightforward modification, almost identical to that in  \cite[Section 3]{Lajda09},
of a standard argument 
from \cite[p.287]{Nev} 
shows that every component $C$
of $G^{-1}(E_N)$ is simply connected, and contains either no zeros of $G$
and one transcendental singularity of $G^{-1}$ over $0$,
or 
exactly one point at which  $G(z) = 0$, which may be a multiple zero.
This is accomplished by deleting from the half-plane ${\rm Re} \, v < \log R$ all pre-images 
of the $\lambda_j$ under $e^v$, and considering $\phi(v) =
G^{-1}(e^v)$ on the resulting simply connected domain $U_0$, the two possible conclusions for $C$ corresponding to whether or not $\phi$ is univalent on $U_0$. 

To prove Proposition \ref{propnumberzeros} 
it now suffices to establish the following lemma. 

\begin{lem}
 \label{lemcountzeros}
If $m \in \{ 0, \ldots, N \}$ and 
$C$ is a component of $G^{-1}(E_m)$ contained in $A$, let $Z_m(C)$ be the number of zeros of $G$ in $C$, counting multiplicities,
plus the number of neighbourhoods  $\Omega_j(0, s, G)$ contained in $C$, and let $Y_m (C)$ be the number of zeros of $G'$ in $C$, again counting multiplicities.
Then 
\begin{equation}
 \label{zeroineq}
Z_m(C) \leq 1 + Y_m(C) .
\end{equation}
\end{lem}
\textit{Proof.}
The lemma will be proved by backwards induction, and (\ref{zeroineq}) clearly holds when $m=N$. 
Now suppose that $0 < m \leq N$, and that (\ref{zeroineq}) holds whenever $C$  is a component of $G^{-1}(E_m)$
contained in $A$.
Let $D$ be a component of $G^{-1}(E_{m-1})$ contained in $A$; the idea of the proof is to delete from $D$ pre-images of $\lambda_m$, thus leaving residual components of $G^{-1}(E_{m})$, to each of which the induction hypothesis can  be applied.
Take all
points $\zeta_j$  in $D$ with $G(\zeta_j) = w_{m}$; each
pre-image of $\lambda_m$ in $D$ contains at least one $\zeta_j$.
If $\zeta_j$ is not a critical point, continuation of $G^{-1}$ along 
$\lambda_{m}$ gives a path $\sigma_j$ from $\zeta_j$ to $\partial D$. These paths $\sigma_j$ are pairwise disjoint, because $G$ has no critical values
on $\lambda_{m}$ apart from $w_{m}$ itself.
Delete all of these $\sigma_j$ from $D$; the set $D'$ which is left is open and is still
connected, because if a path in $D$ joining two points 
of $D'$ meets any of these $\sigma_j$, then it meets only finitely many of them, and may be diverted around each so as not to  leave $D'$.

Next, consider  multiple $w_m$-points of $G$ in $D'$: these are finite in number since $M$ is finite.
Let
  $\zeta_j \in D'$ be a zero of $G - w_{m}$ of multiplicity $m_j + 1 \geq 2$.
Then there are $m_j+1$ pre-images $\tau_{j,k} \subseteq D'$ of $\lambda_{m}$ starting at 
$\zeta_j$ and joining $\zeta_j$ to $\partial D$. Here the $\tau_{j,k}$ for a given $j$ are disjoint, apart from their common starting point
$\zeta_j$, and those starting at distinct $\zeta_j$
do not meet at all. Let $t$ be small and positive and let $T_j = \bigcup_{k=1}^{m_j+1} \tau_{j,k} $; then $U_j = B(\zeta_j, t) \setminus T_j$ has
$m_j+1$ components, and every $ \zeta \in D' \setminus T_j$ can be joined initially to $\zeta_j$ by a path in $D'$,
and hence to a point in $U_j$ by a path in $D' \setminus T_j$. 
It follows that if the $T_j$ are deleted one at a time from $D'$,  each step increases the number of
residual components by at most $m_j$.
Hence the number $r$ of components $C_j$ of $G^{-1}(E_m)$ contained in $D$ exceeds by at most $1$ the number of  zeros of $G'$ 
in $D$ which are also zeros of $G-w_{m}$. It now follows from the induction hypothesis that 
$$
Z_{m-1}(D) \leq \sum_{j=1}^r Z_m(C_j) \leq  \sum_{j=1}^r ( 1 + Y_m(C_j)  )  = r + \sum_{j=1}^r  Y_m(C_j) \leq  1 + Y_{m-1}(D). 
$$
\hfill$\Box$
\vspace{0.1in}

\section{Proof of Theorem \ref{thmUHWV1}}\label{pfuhwv1}

Let the function $L$ be as in the hypotheses of Theorem \ref{thmUHWV1}, and 
assume that 
$L + L'/L$ 
has finitely many non-real zeros. 
The proof follows quite closely the method of \cite[Theorem 1.3]{Lawiman13}.

\begin{lem}
 \label{2013lem5a}
Define 
$F$ by
\begin{equation}
 \label{Lmdef}
F(z) = z - \frac1{L(z)}, \quad
F' = 1 +  \frac{L'}{L^2} .
\end{equation}
Then $F$ is transcendental, $F'$ has finitely many non-real zeros,
and $L$ and $F$ satisfy
\begin{equation}
 \label{2013tsujiest}
\mathfrak{T}(r, L) + \mathfrak{T}(r, F) = O( \log r )
\end{equation}
as $r \to \infty$. Moreover, there exist finitely many $\alpha \in \C $ such that $F(z)$ or $L(z)$
tends to $\alpha$ as $z$ tends to infinity along a path in $\C \setminus \R$.

For real $K > 0$ let
\begin{equation}
H_K = \{ z \in H : |z| > K \} , \quad 
W_K = \{ z \in H : F(z) \in H_K \}.
\label{HKdef}
\end{equation}
Then there exists a large positive real number $K$ such that $F$ has neither critical
nor asymptotic  values in
$H_K $,
and $F$ maps each component of $W_K$
conformally onto $H_K$. 
\end{lem}
\textit{Proof.} First,  $F$ is transcendental because $L$ is.
Since $L$ has finitely many non-real poles, while 
$L + L'/L$ has finitely many non-real zeros, the functions $Q = 1/L$ and $1- Q' = F' = Q(L+L'/L)$ 
have finitely many non-real zeros. Hence 
applying Hayman's alternative \cite[Chapter 3]{Hay2}
to $Q$ as in \cite{BEL},
with the Tsuji  characteristic replacing  that of Nevanlinna, delivers (\ref{2013tsujiest}), 
whereupon Lemma~\ref{directsinglem2} shows that there 
exist finitely many $\alpha \in \C $ such that $F(z)$ or $L(z)$
tends to $\alpha$ as $z$ tends to infinity along a path in $\C \setminus \R$,
and the existence of $K$ follows. 
\hfill$\Box$
\vspace{0.1in}

\begin{lem}\label{lemtheta}
There exist $\theta \in (\pi/4, 3\pi/4 )$ and $N_0 \in \N$ with the following properties.
First, $L$ 
has no critical nor asymptotic
 values in $R^+ = \{ r e^{i \theta } : 0 < r <  + \infty \}$, and none in
$R^- = \{ r e^{-i \theta } : 0 < r < + \infty \}$. 
Next, define $x$ by  $x \sin  \theta = K$, 
with $K$  as  in Lemma \ref{2013lem5a}. 
Then there exist at most $N_0 $ 
points $z$ lying  on the circle $S(0, 2x )$ of centre $0$ and radius $2x$
which satisfy
$L(z) \in R^+$. 
\end{lem}
\textit{Proof.} 
The existence of $\theta$ follows from  Lemma \ref{2013lem5a}. Moreover, if 
$L(z) \in R^+$ for infinitely many $z \in S(0, 2x)$ then  $L_0 = 
e^{-i \theta } L$ satisfies $L_0\left(4x^2/\overline{ \zeta} \right) = \overline{L_0(\zeta)}$,
which is impossible since $L $ is transcendental.
\hfill$\Box$
\vspace{.1in}

\begin{lem}
 \label{2013lemzerosbdy}
Let $D$ be a component of  $W_K$, 
let $a \in \partial D$ be a zero of $L$, and let $\rho$ 
be small and positive. Then $a$ is unique,
and there exists at most one path lying in $D$ and tending to $a$ which is mapped 
by $L$ onto the arc $\Omega_{\theta,\rho} = \{ t e^{i \theta }: \, 0 < t  < \rho \} $.
\end{lem}
\textit{Proof.} 
This is 
 \cite[Lemmas 4.3 and 7.3]{Lawiman13}, but with $L_m$ in the notation of
 \cite{Lawiman13} replaced by $L$, and rests on two facts: first, 
 $F(z) \sim -1/L(z)$ as $z \to a$; second,  since
 $F$ is univalent on $D$,
there is precisely one component of 
$ \left\{ z \in \C : \,  1/2\rho  < |F(z)| < +
\infty , \, \pi/16  < \arg F(z) < 15 \pi /16  \right\}$ in $D$.
\hfill$\Box$
\vspace{0.1in}

The next lemma is \cite[Lemma 7.4]{Lawiman13} and follows from Lemma \ref{2013lem5a}.

\begin{lem}[\cite{Lawiman13}]
 \label{2013lem22}
There exists a positive integer $N_1$ with the following property.
Let $D$ be a component of $W_K$. Then there exist at most $N_1$ components $\Gamma $ of $\partial D$ with $\Gamma \subseteq H$.
\end{lem}
\hfill$\Box$
\vspace{.1in}

The proof of Theorem \ref{thmUHWV1} will now be completed using a combination of ideas from \cite{Larossi,Lawiman13}. Fix a large positive integer $N_2$ and
let $r \in E_1$ be large, where $E_1$ is as in Definition \ref{defnUHWV}. 


\begin{lem}
 \label{secondwvlem}
The set
$Q_r$ in (\ref{Qrdef}) is contained in a component $D$ of $W_K$. 
\end{lem}
\textit{Proof.}  Let $z \in Q_r$. 
Since (\ref{z0def}) and 
(\ref{wv1a})  give
$1/L(z) = o( {\rm Im} \, z )$,
it follows from (\ref{Lmdef}) that 
$|F(z)| > K$ and ${\rm Im} \, F(z) > 0$. 
\hfill$\Box$
\vspace{0.1in}

\begin{lem}
 \label{pointslem}
There exist $S > 0$ and pairwise distinct points 
$w_j$, for $j=1, \ldots, 4N_2$, 
each of large modulus and satisfying $L (w_j) = Se^{i \theta} \in R^+$, where $R^+$ is as in Lemma \ref{lemtheta}, and all lying in the same component $D$ of the set
$W_K$. 
\end{lem}
\textit{Proof.} 
Use
(\ref{wv1a}) to write, on $Q_r$,
$$
\zeta = \log \frac{z}{z_0}, \quad g(z) = \log L(z)  = N(r) \zeta + \log  L(z_0) + o(1).
$$
Since $|N(r) \zeta | = N(r)^{1-\gamma} $ on $ \partial Q_r$, 
Rouch\'e's theorem implies that 
$g(Q_r)$ contains the 
closed disc of centre $ \log L(z_0)    $ and radius $N(r)^{(1-\gamma)/2}$. 
This gives $N_2$ distinct points $w_j \in Q_r$,
all satisfying $L(w_j) =  Se^{i \theta}$ for some large positive $S$, where  $\theta$ is as in Lemma \ref{lemtheta}, and hence $L(w_j) \in R^+$.  
\hfill$\Box$
\vspace{0.1in}

\begin{lem}
 \label{pointslem2}
For $j=1,  \ldots, 4N_2$ choose a component $\sigma_j$ of $L^{-1}(R^+)$ with
$w_j \in \sigma_j$. Then the $\sigma_j$ are pairwise disjoint  and 
each is mapped injectively onto $R^+$ by $L$. 
Moreover at least $2N_2$ of the $\sigma_j$ are such that
$\sigma_j $ lies in $ H_{2x} \cap D$ and  has the following property: as $w \to 0$ on $R^+$ the pre-image 
$z = L^{-1}(w) \in \sigma_j$ tends to infinity in $D$. 
\end{lem}
\textit{Proof.} The first two assertions 
follow from the choice of $\theta$ in Lemma \ref{lemtheta}.
Since the $|w_j|$ and $N_2$ are large,  Lemma \ref{lemtheta} implies that at least $3N_2$ of the 
$\sigma_j$ lie in $H_{2x}$: take $z $ on one of these $\sigma_j$. 
Because $L(z) \in H$, (\ref{Lmdef}) gives $F(z) \in H$. 
If $|L(z)|  \geq 1/x$ then $|F(z)| > 2x - x  > K$
while  $|L(z)| = s < 1/x$ implies that
$$
|F(z)| \geq {\rm Im} \, F(z) \geq \frac{\sin \theta}s > x \sin \theta = K.
$$
Thus at least $3N_2$ of the  $\sigma_j$
lie in $W_K$ and so in $D$. As $L(z) \to 0$ on  one of these
$\sigma_j$, the pre-image  $z$ tends either to infinity or to a zero $a \in \partial D$ of $L$, the 
latter possible for at most one $\sigma_j$ by
Lemma \ref{2013lemzerosbdy}. 
\hfill$\Box$
\vspace{0.1in}

Assume, after re-labelling if necessary,  that  for
$j=1, \ldots, 2N_2$ the path $\sigma_j$ satisfies the conclusions of Lemma \ref{pointslem2},
and let $\sigma_j'$ be the maximal subpath of $\sigma_j$ on which $|L(z)| \leq S$. 
The  $\sigma_j'$ can  be extended to
simple paths $\tau_j$ in $D$, these
pairwise disjoint except for a common starting point $z^* \in D$.
Since $N_2$ is large,  Lemma \ref{2013lem22} gives at least $N_2$ 
pairwise disjoint domains $\Omega_k $, each bounded by two $\tau_j$, and so
by two of the $\sigma_j'$ and a bounded simple path 
$\lambda_k \subseteq D$, such that  the closure of $\Omega_k$ lies in $ D$. 
Because $F(z) \neq \infty$ on $D$, there exists a small positive
$r_k $ such that 
\begin{equation}
\label{omegak}
\hbox{for all $z \in \partial \Omega_k$, either $\arg L(z) = \theta $ or $|L(z)| \geq r_k$.}
\end{equation}

For each $\Omega_k $,  Lemma \ref{2013lem5a} delivers $P_k \in (0, r_k)$
such that the circle $S(0, P_k)$ contains
no critical values of $L$ and no $\alpha \in \C $ such that $L(z) \to \alpha$ 
along a path tending to infinity  in $H$. 
Choose $u_k \in \partial \Omega_k$, lying  on one of the $\sigma_j'$,
with $L(u_k) = P_k e^{i \theta }$,
and continue
$z = L^{-1}(w)$ along $S(0, P_k)$ in the direction taking $z$ into $\Omega_k$.
By (\ref{omegak}) this 
leads to  
$v_k \in \Omega_k$ with $L(v_k) = P_k e^{-i \theta }$.  
Next, Lemma \ref{lemtheta} permits $ L^{-1}(w)$ to be continued along the 
half-ray $w = t e^{-i \theta } $, so that $t$ decreases and  $z = L^{-1}(w)$ starts at $v_k$ and,
 by (\ref{omegak}) again, remains in
$\Omega_k \subseteq D$. 
Since $L(z) \neq 0$ on $D$ this gives a path tending to infinity in 
$\Omega_k$ on which $L(z) \to 0$ with $\arg L(z) = - \theta $. 
Hence 
there exists an unbounded component 
$V_k$ of $\{ z \in \C : {\rm Im } \, (1/L(z)) > 2/P_k \}$, with
$V_k \cup \partial V_k \subseteq \Omega_k \subseteq D$ by (\ref{omegak}) again. Again since $L$ has no zeros in $D$,
the function
$$
u_k(z) = {\rm Im} \, \frac1{L(z)}  -   \frac2{P_k} \quad (z \in V_k), \quad
u_k(z) =  0 \quad (z \not \in V_k) ,
$$
is non-constant and subharmonic in the plane. 
There are at least $N_2$ of these $u_k$,
with disjoint supports, and $N_2$ is large.
Thus the Phragm\'en-Lindel\"of
principle \cite{Hay7} gives, for at least one $k$, a point  $z  \in V_k \subseteq D \subseteq W_K$, with $|z| $ large 
and ${\rm Im} \, 1/L(z)  > |z|^2 $, 
and hence ${\rm Im} \, F(z) < 0$ by (\ref{Lmdef}), which is plainly a contradiction. 
\hfill$\Box$
\vspace{0.1in}


\section{Sufficient conditions for the UHWV property} \label{sufficient}

The main focus of this section will be on proving that if $L$ is as in the 
hypotheses of Theorem \ref{thm1AA} then $L$ has the UHWV property in Definition
\ref{defnUHWV}. 
Let $h(z) = \sum_{n=0}^\infty \alpha_n z^n$ be a transcendental entire function.
Then 
for  $r > 0$ the central index $N(r)$  of $h$  is the largest $n$ for which 
$|\alpha_n| r^n = \max \{ | \alpha_m | r^m : \, m = 0, 1, 2, \ldots \}$, and $N(r)$
tends to infinity with $r$ \cite{Hay5}. The following is a routine consequence
of the Wiman-Valiron theory~\cite{Hay5}. 

\begin{lem}
\label{auxwvlem}
Let  $h$ be a real
transcendental entire function, 
denote by $N(r)$  the central index of $h$ and let 
$1/2 < \tau < \gamma < 1$. 
Then there exists a set $E_0 \subseteq [1, + \infty )$ of finite logarithmic measure 
such that
\begin{equation}
\lim_{r \to + \infty, r \not \in E_0} \frac{N(r)}{\left(  \log M(r, h) \right)^2 } = 0 .
\label{N(r)estwv}
\end{equation}
Furthermore, for each 
$r \not \in  E_0$ there exists $z_0 = z_0(r)$ satisfying (\ref{z0def}), 
such that $|h(z_0)| \sim M(r, h)$ and, uniformly as $r \to + \infty $ outside $E_0$,
\begin{equation}
\label{mainwvest}
h (z) \sim h(z_0) \left( \frac{z}{z_0} \right)^{N(r)} 
\quad \hbox{and} \quad 
\log |h(z)| \geq (1-o(1)) \log M(r, h)
\end{equation}
on the set $Q_r$ in (\ref{Qrdef}). 
\end{lem}
\textit{Proof.}
Choose any $\sigma$ with 
$1/2 < \sigma < \tau < \gamma < 1$.
By a standard result from the Wiman-Valiron theory \cite{Hay5}, 
there  exists a set $E_0 $ of finite logarithmic measure such that (\ref{N(r)estwv}) holds;
furthermore, if
$r \in [1, \infty) \setminus E_0$ and $|z_1| = r$, $|h(z_1)| \sim M(r, h)= \max \{ |h(z)| : |z| = r \} $,
then
$$
h (z) \sim h(z_1) \left( \frac{z}{z_1} \right)^{N(r)}
\quad \hbox{for} \quad \left| \log \frac{z}{z_1} \right| \leq N(r)^{-\sigma} .
$$
Since $h$ is real, 
it may  be assumed that ${\rm Im} \, z_1 \geq 0$ for $r \in [1, \infty) \setminus E_0$,
and so 
there exists $z_0$ satisfying (\ref{z0def}) and the first estimate of
(\ref{mainwvest}).  Next, (\ref{N(r)estwv}) and 
the fact that $0 < 1 - \gamma < 1/2$ yield
 the second estimate of (\ref{mainwvest}) via
 \begin{equation*}
\label{hminest}
\log |h(z)| \geq  \log M(r, h) -   N(r)^{1-\gamma}  - o(1) \geq (1-o(1)) \log M(r, h) .
\end{equation*} 
\hfill$\Box$
\vspace{0.1in}

Lemma \ref{auxwvlem} shows that every real transcendental entire function has the UHWV property.
The same is in fact true of any real transcendental meromorphic function in the plane
for which the inverse function has a direct transcendental singularity over $\infty$: this 
can be proved 
identically, but using
the version of Wiman-Valiron theory developed in \cite{BRS} for functions with direct tracts.

\begin{lem}
\label{mainwvlem}
Let 
$L$ be as in the hypotheses of Theorem \ref{thm1AA}.
Denote by $N(r)$  the central index of $h$ and let 
$1/2 < \tau < \gamma < 1$. 
Then there exists a set $E_0 \subseteq [1, + \infty )$ of finite logarithmic measure 
with the following property: for each 
$r \in  [1, + \infty) \setminus E_0$ there exists $z_0 = z_0(r)$ satisfying (\ref{z0def}) 
such that (\ref{wv1a}) holds on the set $Q_r$ in (\ref{Qrdef}), uniformly as $r \to + \infty $ outside $E_0$.
In particular, $L$ has the UHWV property, with $E_1 = [1, + \infty ) \setminus E_0$. 
\end{lem}
\textit{Proof.}
Choose $E_0$ and $z_0$ as in Lemma \ref{auxwvlem}. 
Combining  (\ref{C1}), (\ref{N(r)estwv}) and (\ref{mainwvest}) 
shows that, for large $r \in [1, \infty) \setminus E_0$ and $z \in Q_r$,
\begin{eqnarray*}
R_1(z) \psi_1(z) &\sim& R_1(z_0) \psi_1( z_0), \\
\log \frac1{ |R_1(z) \psi_1(z)| } + \log^+ |R_2(z) \psi_2(z)|
 &\leq& O( \log r) + O( \log N(r) ) \leq   o( \log M(r, h) ), \\
\log |L(z)| &\geq& (1-o(1)) \log M(r, h) ,
\end{eqnarray*}
from which (\ref{wv1a}) follows. 
\hfill$\Box$
\vspace{0.1in}

\begin{lem}
\label{lem1Aimplies1AA}
Let 
$g = g_1/g_2 $,
where $g_1, g_2$ are real entire functions, with no common zeros and finitely many non-real zeros. Assume 
that  $g_2$ has finite lower
order, but $g$ has infinite order. 
Then $L = g'/g$  satisfies the hypotheses of Theorem \ref{thm1AA}, and hence the
conclusions of Lemma \ref{mainwvlem}.
\end{lem}
\textit{Proof.}
The logarithmic derivative of each $g_j$ has a Levin-Ostrovskii factorisation 
\cite{BEL,LeO}
\begin{equation}
\frac{g_j'}{g_j} = \phi_j \psi_j ,
\label{LeOrepgj}
\end{equation}
in which: $\phi_j$ and $\psi_j$ are real meromorphic functions; $\phi_j$ has 
finitely many poles;   if $g_j$ has finitely many zeros  then $\psi_j \equiv 1$;
if $g_j$ has
infinitely many zeros then (\ref{LeOrepgj}) is obtained by
applying  Lemma~\ref{lemlevost} in conjunction
with Rolle's theorem, in which case $\psi_j(H) \subseteq H$. 
Since $g_2$ has finite lower order, (\ref{C1}) and the lemma of the logarithmic 
derivative give
$$
T(r, \phi_2) \leq 
m(r, \phi_2) + O( \log r) \leq m(r, 1/\psi_2) + m(r, g_2'/g_2) + O( \log r) = O( \log r ) 
$$
on a sequence of $r$ tending to infinity, and so $\phi_2$ is a rational function. 
Thus \cite[Lemma 5.1]{BEL} implies that $g_2$ has finite order. Because $g$ has
infinite order, so has $g_1$, and applying  \cite[Lemma 5.1]{BEL} again shows that
$\phi_1 $ is transcendental. Hence $L $ has a representation (\ref{Lnewrep})
as required. 
\hfill$\Box$
\vspace{0.1in}

\subsection{Proof of Theorems \ref{thm1A} and \ref{thm1AA} } \label{pfthm1AA}

First, let $L$ be as in the hypotheses of Theorem \ref{thm1AA}: then 
$L$ has the UHWV property, by Lemma~\ref{mainwvlem}, whereupon Theorem \ref{thmUHWV1}
implies that  $L+L'/L$ has infinitely
many non-real zeros.
Next, if $f$ is as in the hypotheses of  Theorem \ref{thm1A} then $L = f'/f$ satisfies
those of 
Theorem \ref{thm1AA}, by  Lemma \ref{lem1Aimplies1AA} applied to $f$ or $1/f$. 
\hfill$\Box$
\vspace{0.1in}


\section{Proof of Theorem \ref{thmlindiffpol}}\label{pfthmlindiffpol}

Let $L, \phi, \psi, a, b$ be as in the hypotheses of Theorem \ref{thmlindiffpol} and 
suppose that $L' + a L^2 +  b $ has finitely many non-real zeros.
Writing $L(z) = \alpha L_1( \beta z)$, where $\alpha = \sqrt{b/a}$ and $\beta = \sqrt{ab}$,
makes it possible to assume that $a=b=1$.

The following estimate for the Tsuji characteristic of $L$ was deduced 
in \cite[Lemma 4.1]{Lajda09} from an argument of Tumura-Clunie type \cite[Ch. 3]{Hay2}. Note that  \cite[Lemma 4.1]{Lajda09}
is stated for the special case in which $L = f'/f$, where $f$ is an entire function such that $f$ and $f''+f$ have finitely many non-real zeros, 
but the proof depends only on $L$ having finitely many non-real poles and  $L'+L^2+1$   finitely many non-real zeros. 

\begin{lem}
[Lemma 4.1, \cite{Lajda09}]
\label{lemA1}
The Tsuji characteristic of  $L$ satisfies 
\begin{equation}
\mathfrak{T}(r, L) = O( \log r ) \quad \hbox{as $r \to \infty$}.
\label{A2}
\end{equation}
\end{lem}
\hfill$\Box$
\vspace{0.1in}

\begin{lem}
 \label{phitrans}
The transcendental entire function $h$  in (\ref{Lnewrep}) has order at most $1$. 
\end{lem}
\textit{Proof.} This is fairly standard.
First, (\ref{Lnewrep}) and (\ref{C1}) imply that, as $r \to + \infty$, 
$$
T(r, h) = m(r, h) \leq m(r, L) + O( \log  r ) .
$$
This implies in turn that, as $R \to + \infty $, 
by (\ref{A2}) and 
an inequality of Levin and Ostrovskii \cite[p.332]{LeO}
(see also \cite[Lemma~3.2]{BEL} or \cite[Lemma 2.3]{Lajda09}), 
$$
\frac{T(R, h)}{2R^2} \leq 
\int_R^\infty \frac{T(r, h)}{r^3} \, dr \leq 2
\int_R^\infty \frac{\mathfrak{T}(r, L)}{r^2 } \, dr + 
O \left( \frac{\log R}{R^2} \right) = O \left( \frac{\log R}{R} \right) .
$$
\hfill$\Box$
\vspace{0.1in}

The proof in \cite{Lajda09} made extensive use of the auxiliary function 
$F = (TL-1)/(L+T)$, where $T(z) = \tan z $.
For the present paper it turns out to be simpler to work with
\begin{equation}
 \label{g1}
G(z) = e^{2iz} \left( \frac{L(z)-i}{L(z)+i} \right) = - \left( \frac{ F(z)-i}{F(z)+i}  \right) ,
\quad    G'(z) = \frac{ 2i e^{2iz} (L'(z)+L(z)^2 + 1)}{(L(z)+i)^2} .
\end{equation}
Then $|G(x)| = 1$ for $x \in \R$, and $G'$ 
has finitely many zeros in $\C \setminus \R$, while
\begin{equation}
\label{YWdef}
 Y = \{ z \in H: \, L(z) \in H \} \subseteq W = \{ z \in H: \, |G(z)| < 1 \}. 
\end{equation}

There now follows a sequence of lemmas which together show that $G$ has finitely many  asymptotic values $\alpha \in \C $ 
with $| \alpha | \neq  1$, 
using a method which substantially simplifies the approach in \cite{Lajda09}. 
For $\alpha \in \C$, use (\ref{g1}) to define $s_\alpha$ by 
\begin{equation}
 \label{salphadef}
s_\alpha(z) = 
\frac{ G(z) - \alpha}{e^{2iz} - G(z)} 
= \frac1{2i}  \left( L(z)-i  - \alpha e^{-2iz} (L(z)+i) \right)  .
\end{equation}
Since $L$ has finitely many non-real poles, so has each  $s_\alpha$.
 
\begin{lem}
\label{lems1}
 Let $\alpha, \beta \in \C$ satisfy $\alpha \neq \beta$.
Then there exists $c_1 > 0$ such that if $z \in H$ and $|z|$ is large 
then $| s_\alpha(z) | + | s_\beta(z) | \geq  c_1$.
\end{lem}
\textit{Proof.} 
Assume that there exists a sequence $z_n \to \infty$ in $H$ such that 
$| s_\alpha(z_n) | + | s_\beta(z_n) | \to 0$. Since $|e^{2i z_n} | \leq 1$ in (\ref{salphadef}), it must be the case that $G(z_n) = O(1)$, 
from which it follows that $G(z_n) \to \alpha$ and $G(z_n) \to \beta$,
which is impossible. 
\hfill$\Box$
\vspace{.1in}

\begin{lem}
\label{lems2}
 Let $\alpha, \beta \in \C$ satisfy $\alpha \neq \beta$, and let $c_2 > 0$. 
Then there exists $c_3 > 0$ such that if $z_n \to \infty$ in $H$ with $|e^{2iz_n} - \alpha | \geq c_2$ 
and $G(z_n) \to \alpha$ 
then $s_\alpha (z_n) \to 0$ and $| s_\beta(z_n)| \leq c_3$. 
\end{lem}
\textit{Proof.} This follows  from (\ref{salphadef}) and the fact that $2 |e^{2iz_n} - G(z_n)| \geq c_2$ for all large $n$. 
\hfill$\Box$
\vspace{.1in}

\begin{lem}
\label{lems3}
 Let $\alpha_1, \dots, \alpha_N  \in \C$
be pairwise distinct, with $| \alpha_1 | \neq 0, 1$, and let 
$G(z) \to \alpha_1$ on a path $\gamma$ tending to infinity in $H$. 
Then there exists a path $\lambda $  tending to infinity in $H$ on which 
$s_{\alpha_1}(z) \to 0$ and $s_{\alpha_j}(z)$ is bounded for $j=2, \ldots, N$.
\end{lem}
\textit{Proof.} Evidently there exists $q \in \C \setminus \R$ such that  the solutions of $e^{2iz} = \alpha_1 $ are
$a_n = n \pi + q $, $n \in \Z$. Let $\varepsilon$ be small and positive. 
Then Lemma \ref{lems2} shows that $s_{\alpha_1}(z)$ is small, and the remaining $s_{\alpha_j}(z)$ are uniformly bounded, for all $z \in \gamma$ 
such that $|z|$ is large and $z$ lies outside the union of the  discs $B(a_n, \varepsilon)$.

It may therefore be assumed that 
$\gamma$ meets the  disc $B(a_n, \varepsilon)$ for all $n$ in an unbounded set $E \subseteq \Z$, since otherwise there is nothing further
to prove.
Then $0 < | \alpha_1 | < 1$ and for each $n \in E$
there exists a simple subpath $\sigma_n$ 
of $\gamma$ which lies in the annulus $ 2 \varepsilon  \leq |z-a_n| \leq  4 \varepsilon$ and
joins the two boundary circles.
Lemma \ref{lems2} implies that 
\begin{equation}
 \label{sigmafact}
\lim_{|n| \to \infty, n \in E}  \tau_n = 0,
\quad \tau_n = \max \{ |G(z) - \alpha_1 | + |s_{\alpha_1}(z)| : \, z \in \sigma_n \} .
\end{equation}
Moreover, standard estimates \cite{Nev} give a positive  constant 
$C$, independent of $n \in E$,  such that the harmonic measure $\omega(z, \sigma_n, B(a_n, 4\varepsilon) \setminus \sigma_n)$ is at least $C$
for $|z-a_n| \leq \varepsilon$.

Let $E_1$ be the set of $n \in E$ such that  $|n|$ is large and there exists  $z_1$ in $B(a_n, 4 \varepsilon)$ with $|L(z_1)| \leq 2$. 
Since the functions $L_n(z) = L(a_n + z)$, $n \in E_1$, satisfy $L_n \neq \infty $ and $L_n'+L_n^2+1 \neq 0$ on $B(0, 8 \varepsilon)$, 
Lemma \ref{thmel} and  (\ref{salphadef})
deliver  $K_1, K_2  > 0$, independent of $n$, such that 
$|L(z)| \leq K_1$ and $| s_{\alpha_j}(z)| \leq K_2$
for $z$ in $B(a_n, 4 \varepsilon)$, $n \in E_1$  and  $j=1, \ldots, N$.
This makes $u_1(z) =  \log |s_{\alpha_1}(z)/K_2|$ subharmonic and non-positive on $B(a_n, 4 \varepsilon)$, for $n \in E_1$, and a standard
combination of (\ref{sigmafact}) with the two constants theorem \cite{Nev} yields,  for $|z-a_n| \leq \varepsilon$,
$$
u_1(z) \leq C \log \left( \frac{\tau_n}{K_2} \right) , \quad |s_{\alpha_1}(z)| \leq K_2 \left( \frac{\tau_n}{K_2} \right)^C .
$$
Thus for  $n \in E_1$ and $z \in \gamma \cap B(a_n, \varepsilon)$, the term  $s_{\alpha_1}(z)$ is small,  by (\ref{sigmafact}), 
while $| s_{\alpha_j}(z)| \leq K_2$ for $j=2, \ldots, N$.

It remains only to deal with the set $E_2$ of $n \in E \setminus E_1$ such that $|n|$ is large. These $n$ are such that 
$|L(z)| > 2$ for all $z$ in $B(a_n, 4 \varepsilon)$, and hence $|G(z)| \leq 3$ there, by (\ref{g1}).
This time $u_2(z) = \log |(G(z)- \alpha_1)/4| $ is subharmonic and non-positive on $B(a_n, 4 \varepsilon)$, and 
combining (\ref{sigmafact}) with the two constants theorem  yields
$\displaystyle{|G(z)- \alpha_1| \leq 4 \left( \frac{\tau_n}{4} \right)^C}$  for $|z-a_n| \leq \varepsilon$.
Thus for  $|z-a_n| = \varepsilon$, where $n \in E \setminus E_1$ and $|n|$ is large,   (\ref{sigmafact}) and Lemma \ref{lems2} imply that  $s_{\alpha_1}(z)$ 
is small, and the remaining $s_{\alpha_j}(z)$ are uniformly bounded. The proof is now completed by replacing 
any part of $\gamma$ which enters and leaves $B(a_n, \varepsilon)$, for  $n \in  E_2$, 
by an arc of the circle $| z- a_n| =  \varepsilon$.
\hfill$\Box$
\vspace{.1in}

\begin{lem}
\label{lems4}
The function $G$ has finitely many  asymptotic values $\alpha \in \C $ 
with $| \alpha | \neq 1$. 
\end{lem}
\textit{Proof.} Since $|G(x)| = 1$ for $x \in \R$ it suffices to show that there do not exist pairwise
distinct $\alpha_1, \alpha_2, \alpha_3  \in \C$,
with   $| \alpha_j | \neq 0, 1$, 
such that $G(z) \to \alpha_j$ along a path $\gamma_j$ tending to infinity in $H$. 
Assume the contrary: then Lemma \ref{lems3} gives paths $\lambda_1, \lambda_2, \lambda_3$ in $H$ such that
$s_{\alpha_j}(z)$ tends to $0$, while the remaining $s_{\alpha_k}(z)$ are bounded,   as $z \to \infty$ on $\lambda_j$.  Hence 
$Q(z) = s_{\alpha_1}(z) s_{\alpha_{2}}(z) s_{\alpha_{3}}(z)$ tends to $0$  on each $\lambda_j$. By Lemma \ref{lems1}, each intersection $\lambda_j \cap \lambda_{j'}$
is bounded for $j \neq j'$. 

Choose a large $R \in (0, \infty)$. It may be assumed that the $\lambda_j$ start on $|z| = R$ and divide $\{ z \in H : \, |z| > R \}$ into four 
disjoint unbounded domains $D_0, \ldots, D_3$, such that $\lambda_j$ separates  $D_{j-1}$ from $D_j$ for $j =1, 2, 3$. 
Suppose first that,  as $z$ tends to infinity in $D_1$, the function  $Q(z)  $ 
is bounded, and so tends to $0$
by  the Phragm\'en-Lindel\"of principle.  Lemma \ref{lems1} implies that $|s_{\alpha_2}(z) s_{\alpha_3}(z)| >   |s_{\alpha_1}(z)|     $ 
as $z \to \infty$ on $\lambda_1$, while $|s_{\alpha_1}(z)| >    |s_{\alpha_2}(z) s_{\alpha_3}(z)|     $ 
as $z \to \infty$ on $\lambda_2$. Hence there exists $z \in D_1$, with $|z|$ arbitrarily large, $Q(z)$ small
and $|s_{\alpha_1}(z)| = |s_{\alpha_2}(z) s_{\alpha_3}(z)|$. But this implies that $s_{\alpha_1}(z)$ and at least one of
$s_{\alpha_2}(z), \,  s_{\alpha_3}(z)$ are both small, which contradicts Lemma \ref{lems1}. 

It follows that $Q(z)$ is unbounded   as $z$ tends to infinity in $D_1$ and, by the same argument, in $D_2$ also, so that $Q^{-1}$ has
at least two direct singularities over $\infty$, lying in $H$. Since $\mathfrak{T}(r, Q) = O( \log r )$ as $r \to \infty$,
by (\ref{A2}) and (\ref{salphadef}), this  contradicts Lemma \ref{directsinglem}. 
\hfill$\Box$
\vspace{.1in}

\begin{lem}
\label{lems5}
If $a \in \C$ and $|a| \neq 1$, and if $G^{-1}$ has a transcendental singularity over $a $,
then the singularity is direct.
\end{lem}
\textit{Proof.} This follows from Lemma \ref{lems4}, the fact that $G'$ has finitely many non-real zeros, and the 
standard classification of isolated singularities of the inverse function \cite[p.287]{Nev}. 
\hfill$\Box$
\vspace{.1in}


\begin{lem}\label{lemB1}
If  $G^{-1}$ has a  transcendental singularity  
over $a \in \C$,
then $a = 0$ or $|a|= 1$. Moreover, there exists $N_0 \in \N$ such that if $D$ is a component of 
the set $W$ in (\ref{YWdef}) then $G$ has at most $N_0$ zeros in $D$, counting multiplicities.

Furthermore,  $L$ has finitely many non-real asymptotic values and $L^{-1}$
cannot have a direct transcendental singularity 
over $a \in \C \setminus \R$. 
Finally, there exists  $\theta \in (\pi/4, 3 \pi/4)$ such that $L$ has no critical or asymptotic values on the open half-line $R^+$ given by
$w = i + t e^{i \theta} $, $0 < t < + \infty$. 
\end{lem}
\textit{Proof.} 
This is  a modification of  \cite[Lemma 5.1]{Lajda09}. 
Assume first that $g$ is $G$ or $L$ and  that
 $g^{-1}$ has a direct transcendental singularity 
over $a \in \C$, with 
$|a| \neq 0,  1$ if $g = G$, and $a \in \C \setminus \R$ if $g = L$. 
Since $|G(x)| = 1$ on $\R$ and $L(\R) \subseteq \R \cup \{ \infty \}$, 
it may be assumed that the singularity lies in $H$.
Let $\delta_1, \delta_2 $ be small and positive.
Then there
exists a component $D \subseteq H$ of 
$\{ z \in \C : |g(z) - a| < \delta_1 \}$ such that $g(z) \neq a$ on $D$ and 
$$
v(z) = \log \frac{\delta_1}{|g(z) - a|}  \quad (z \in D), \quad
v(z) = 0 \quad (z \in \C \setminus D),
$$
defines a non-constant subharmonic function on $\C$. 
Because $\mathfrak{T}(r, g) = O( \log r )$
as $r \to \infty$ by (\ref{A2}) and (\ref{g1}), a standard argument as in  \cite[(2.2)]{Lajda09} shows that
$v$ has order of growth 
at most $1$. 

Let $N_1$ be a large positive integer. 
By Lemma \ref{phitrans},
there exists a real polynomial $P_1$, of degree at most $N_1-1$, 
such that  $h_1(z) = z^{-N_1}( h(z) - P_1(z))$
is entire and transcendental of
order at most $1$. Let $C$ be a component of the set $\{ z \in \C : |h_1(z)| > 1 \}$.
If $z \in C$ and $|z|$ is large, and $\delta_2 <  | \arg z |  < \pi - \delta_2$,
then  combining (\ref{Lnewrep}), (\ref{C1}) and (\ref{g1}) shows that
$|L(z)| = |L( \bar z)|$ is large, while  one of $|G(z)|$ and $|G(\bar z)|$  is small and
the other is large. 
Thus neither $z$ nor $\bar z$ lies in $D$,
and so $z$ cannot lie in the reflection of $D$ across $\R$.
For $s > 0$ denote by $\theta_C (s), \theta_D (s)$  the
angular measure of the intersection of $C$, respectively $D$, with the circle $|z| = s$,
and let $\theta_C^*(s) = + \infty$ if the whole circle $|z| = s$ lies in $C$, with 
$\theta_C^*(s) = \theta_C(s) $ otherwise. 
Then the Cauchy-Schwarz inequality and the fact that $\delta_2$ is small
give, for all large~$s$,
$$\theta_C(s) + 2 \theta_D(s) \leq 2 \pi + 8 \delta_2 ,
\quad 9 \leq 
\left( \frac1{ \theta_C^*(s)} + \frac2 {\theta_D(s)} \right) ( 2 \pi + 8 \delta_2 ). $$
Now $v_1 = \log^+ |h_1|$ and $v_2 = v$ are both subharmonic 
of order at most $1$, so set $B^*(r, v_j) =  \max \{ v_j(z) : \, |z| = r \}$ for $r > 0$,
 fix a large positive $r_0$ and let $r \to + \infty$. 
Then
a standard application of  Carleman's estimate for harmonic measure \cite{Tsuji}, exactly
as in \cite[Lemmas 2.1 and 5.1]{Lajda09}, 
leads to a contradiction via
\begin{eqnarray*}
\left( \frac{9}{2  + 8 \delta_2 /\pi } \right) \log r &\leq & 
\int_{r_0}^r \left( \frac\pi{s  \theta_C^*(s)} + \frac{2\pi} {s \theta_D(s)} \right) \, ds + O(1) \\
&\leq& 
\log B^*(2r, v_1)  + 2 
\log B^*(2r, v_2)  + O(1) \\
&\leq& (3+o(1)) \log r . 
 \end{eqnarray*} 

In view of Lemma \ref{lems5}, this shows
 that if  $a \in \C$  is an asymptotic value of $G$ then $a = 0$ or $|a|= 1$. Hence the integer  $N_0$ exists by
Proposition \ref{propnumberzeros} and  the fact that 
$G'$ has finitely many non-real zeros.
Next,  $L^{-1}$
cannot have a direct transcendental singularity 
over $a \in \C \setminus \R$, and
$L$ cannot have infinitely many non-real asymptotic values, by  (\ref{A2}) and Lemma~\ref{directsinglem2}
(indeed,   
$L^{-1}$ would otherwise have at least two direct transcendental singularities
over $\infty$, lying in $H$, contradicting
Lemma~\ref{directsinglem}).
The existence of $\theta$ follows at once. 
\hfill$\Box$
\vspace{0.1in}

Since $L$ satisfies the hypotheses of Theorem \ref{thmlindiffpol}, which are the same as those of Theorem~\ref{thm1AA},  Lemma~\ref{mainwvlem} may now be applied to $L$.

\begin{lem}
 \label{secondwvlema}
 Let $r \not \in E_0$ be large, with $E_0$ as in Lemma~\ref{mainwvlem}. Then
the set  $Q_r$ in (\ref{Qrdef}) is contained in a component $D$ of the set $W$ in
(\ref{YWdef}). 
\end{lem}
\textit{Proof.} Since $r$ is large, (\ref{z0def}), (\ref{wv1a}) and (\ref{g1}) imply that, for $z \in Q_r$,
\begin{eqnarray*}
{\rm Im } \, z &\geq&   N(r)^{-\tau }, \quad \frac1{L(z)} = o \left( N(r)^{-\tau } \right),\\
 G(z) &=& e^{2iz} \cdot \frac{1 -i/L(z)}{1+i/L(z)}  = e^{2iz} ( 1 + \varepsilon_1(z)), \quad 
 \varepsilon_1(z) = o \left( N(r)^{-\tau } \right),\\
\log |G(z)| &\leq& - 2 N(r)^{-\tau} +  o \left( N(r)^{-\tau } \right) < 0 .
\end{eqnarray*}
\hfill$\Box$
\vspace{0.1in}

\begin{lem}
 \label{pointslemA}
Let $N_2$ be a large positive integer. Then for large enough $r$  as in Lemma \ref{secondwvlema} there exist $S > 0$ and  pairwise distinct 
$w_j \in Q_r$, for $j=1, \ldots, N_2$, 
such that  $L(w_j) = i + S e^{i \theta } \in R^+$, where $\theta$ and $R^+$ are as in Lemma \ref{lemB1}. 

For each $j \in \{ 1, \ldots, N_2 \}$, let  $\sigma_j$ be the  component of $L^{-1}(R^+)$ with
$w_j \in \sigma_j$. Then the $\sigma_j$ are pairwise disjoint  and lie in the same component $D$ of $W$ as $Q_r$, and
each is mapped injectively onto $R^+$ by $L$. 
Furthermore, at least one of the $\sigma_j$ has the property that
as $w \to i$ on $R^+$ the pre-image 
$z = L^{-1}(w) \in \sigma_j$ tends to infinity in $D$. 
\end{lem}
\textit{Proof.} Let $r$ be large and as in Lemma \ref{secondwvlema}. 
The existence of $S$ and the $w_j$ is proved exactly as in Lemma \ref{pointslem}, 
using the fact that $L(z)-i \sim L(z)$ on $Q_r$. 
The next three assertions follow from the fact that $L^{-1}$ has no singular values on $R^+$, by the
choice of $\theta$, and the inclusions $w_j \in Q_r \subseteq D$ and (\ref{YWdef}). Now, as $w \to i$ on $R^+$ the pre-image 
$z = L^{-1}(w) \in \sigma_j$ lies in $D$ and tends either  to  a zero of $L-i$, which by (\ref{g1}) is a zero of $G$ in $D$ of the same multiplicity,
or to infinity. 
Because $N_2 $ is large, Lemma \ref{lemB1} now implies that 
$z= L^{-1}(w)$ must tend to infinity for at least one  $j$. 
\hfill$\Box$
\vspace{0.1in}

The proof of Theorem~\ref{thmlindiffpol} may now be completed. 
Lemma 
\ref{pointslemA} shows that
$L(z) $ tends to~$i$ along a path $\mu$ tending to infinity in the component $D$ of $W$.
This gives $t \in (0, 1/2)$ and a neighbourhood $\Omega(t)$ of a transcendental singularity of $L^{-1}$ over $i$,
such that $\mu \setminus \Omega(t)$ is bounded. Moreover,
$\Omega(t)$ lies in $H$, and so in $Y \subseteq W$, by (\ref{YWdef}), and hence in $D$.
By Lemma~\ref{lemB1}  and  (\ref{g1}),
$G$ and $L-i$ have finitely many zeros in $D$. But this implies that $L^{-1}$ has a direct  transcendental singularity  over $i$, 
which contradicts Lemma \ref{lemB1}.
\hfill$\Box$
\vspace{0.1in}

\section{Proof of Theorem \ref{cor1}}

Assume that $f$ satisfies the hypotheses of Theorem \ref{cor1}, and hence those of
Theorem \ref{thm1A}. Then, as in Section \ref{pfthm1AA},
$L = f'/f$ satisfies the identical hypotheses
of Theorems \ref{thm1AA} and  \ref{thmlindiffpol}, and the latter 
 implies that
$L'+L^2 + \omega = (f''+\omega f)/f$ has infinitely many non-real zeros. 
\hfill$\Box$
\vspace{0.1in}

\section{Proof of Theorem
 \ref{thm11}} 
 \label{infiniteorder}

Let $U$ be as in the hypotheses of Theorem \ref{thm11}. 
Part (i) was already proved in Lemma \ref{lemlevost},
 with $\psi$ as in (\ref{wwpsi}).
Next,
assume that $S$ has infinite order in (\ref{ww1}) and write 
$$
\frac{U'}U = L +  \frac{\psi'}{\psi} ,  \quad  L = \frac{S'}S .
$$
Since $S$ has finitely many poles and non-real zeros, 
applying Lemma \ref{lem1Aimplies1AA}, with $g = S$,  
shows that $L$ satisfies the hypotheses of Theorem \ref{thm1AA},
and therefore has the UHWV property by Lemma \ref{mainwvlem}.
Then the UHWV property for $U'/U$
follows 
from  the fact that (\ref{z0def}), (\ref{wv1a}) and (\ref{C1})   give, 
for large $r \in E_1$ and $z \in Q_r$, 
$$
\left| \frac{\psi'(z)}{\psi(z)} \right| = O\left(  \frac{N(r)^\tau }{r} \right) = o( |L(z)| ), \quad 
\frac{U'(z)}{U(z)} \sim L(z)  . $$
Since  $U'/U$ has finitely many non-real poles,
part (ii) of Theorem \ref{thm11} follows from
Theorem~\ref{thmUHWV1}.


  Assume henceforth that
  $S = Re^P = U/\psi$ is  as in the hypotheses of part (iii), 
  in particular with $P$ a non-constant polynomial, and let $2 \leq m \in \N$.
  Since $\psi$ maps $H$ into itself,  (\ref{wwpsi})  
  gives a series representation
 \begin{equation}
 \psi(z) = Az + B + \sum_{k \in K}  A_k \left( \frac1{x_k - z} - \frac1{x_k} \right) ,
 \label{ww3}
 \end{equation}
 with $A, B, A_k \in \R$,  $A \geq 0$, $A_k > 0$ and 
 $\sum_{k \in K} A_k x_k^{-2} < + \infty $ \cite{Le}.  
 Write, 
 using~(\ref{ww1}),
 \begin{equation}
 U^{(m)} = S \sum_{j=0}^m a_j \psi^{(j)} = S a_0 \phi ,
 \quad 
 a_j =  { m \choose j } \frac{S^{(m-j)}}S = a_0 b_j,
 \quad \phi =  \phi_m =  \sum_{j=0}^m b_j \psi^{(j)} .
  \label{ww4}
 \end{equation}
 Here $a_m = b_0 = 1$ and, by  
 standard estimates based on the formulas
 $$
 \frac{S'(z)}{S(z)} = \frac{R'(z)}{R(z)} + P'(z), \quad 
 \frac{S^{(p+1)}(z)}{S(z)} = \frac{S^{(p)}(z)}{S(z)} \cdot \frac{S'(z)}{S(z)}  + \frac{d}{dz}
 \left(   \frac{S^{(p)}(z)}{S(z)} \right), 
 $$
 the real rational functions $a_j, b_j$ satisfy 
 \begin{equation}
 a_j(z) = { m \choose j }  P'(z)^{m-j}  \left( 1 + O \left( \frac1{|z|} \right) \right) , \quad 
  b_j(z) = { m \choose j }  P'(z)^{-j}  \left( 1 + O \left( \frac1{|z|} \right) \right) 
 \label{ww7}
 \end{equation}
as $z \to \infty$.  
 The key to the proof of Theorem \ref{thm11}(iii) is the following. 
 
 \begin{prop}
 \label{prop1}
 Let  $s_0$ be a large positive real number and let
 $I \subseteq \R  \setminus [-s_0, s_0]$ be an open interval  which contains no
 poles of $\psi$. Then the number of zeros of $U^{(m)}$ in $I$, 
 counting multiplicities, is at most $1$ if $m$ is even, and at most $2$ if $m$ is odd. 
 
 Next, let $k \in K$ be such that $|k|$ is large, and let 
$n_{k,m}$ be the number of zeros of
$U^{(m)}$, counting multiplicities, in $(x_k, x_{k+1})$. 
If $m$ is even then $n_{k,m} = 1$, while if $m$ is odd then $n_{k,m} \in \{ 0, 2 \}$.
 \end{prop}
 \textit{Proof.} 
 Suppose first that $m$ is even,
 set $b_{m+1} = 0$ and recall that $b_0 = 1$. 
 Since $P$ is a real polynomial, $S$ and $a_0$ do not change sign on $I$ and 
 so (\ref{ww4}) implies that $U^{(m)}$ has the same number of zeros in $I$ as $\phi$.
Thus to prove the first assertion it suffices to show that if 
$s_0$ is sufficiently large then the derivative $\phi'$ is positive on $I$,
where $\phi'$ is given by 
\begin{equation}
\label{ww9a}
\phi ' =  \sum_{j=0}^m \left( b_j \psi^{(j+1)} + b_j' \psi^{(j)} \right) =
 \sum_{j=0}^m  b_j \psi^{(j+1)} + \sum_{j=1}^{m+1} b_{j}' \psi^{(j)} =
 \sum_{j=0}^{m} c_j \psi^{(j+1)} ,
\end{equation}
in which the $c_j$ satisfy, as $z \to \infty$, by (\ref{ww7}), 
\begin{equation}
c_j (z) = b_{j}(z)  + b_{j+1}'(z) 
= { m \choose j } P'(z)^{-j}  \left( 1 + O \left( \frac1{|z|} \right) \right) .
 \label{ww9}
 \end{equation}
 For $x \in I$ let  $X_k = P'(x) (x_k-x) \in \R \setminus \{ 0 \}$
and let $Q_m$ be as in Lemma 
\ref{lem2}. Then  (\ref{ww3}), (\ref{ww9a}), 
(\ref{ww9}) and the fact that $c_0 (\infty) = 1$
deliver the following, in which the $o(1)$ terms are
uniformly small for $x \in I$, provided $s_0$ is large enough:  
\begin{eqnarray*}
\phi'(x) 
&=& \sum_{j=0}^{m} c_j (x) \frac{d^{j+1}}{dx^{j+1}}
 \left (  Ax + B + \sum_{k \in K} A_k \left( \frac1{x_k - x} - \frac1{x_k} \right) \right) \\
 &\geq& \sum_{j=0}^{m} c_j (x) \frac{d^{j+1}}{dx^{j+1}}
 \left (   \sum_{k \in K} A_k \left( \frac1{x_k - x} - \frac1{x_k} \right) \right) \\
   &=&
   \sum_{k \in K} A_k \sum_{j=0}^{m} c_j (x) \frac{d^{j+1}}{dx^{j+1}}
  \left (    \frac1{x_k - x}  \right)\\
  &=&
   \sum_{k \in K} A_k \sum_{j=0}^{m} { m \choose j } P'(x)^{-j}  (1+o(1)) 
  \left (    \frac{(j+1)!}{(x_k - x)^{j+2}}  \right)  \\
  &=&
    \sum_{k \in K} \frac{A_k}{ (x_k-x)^{2} }
    \sum_{j=0}^{m}  (1+o(1))  { m \choose j } \frac{(j+1)! }{X_k^{j} } \\
&=&
\sum_{k \in K} \frac{A_k }{(x_k-x)^{2}} \left( Q_m\left(\frac1{X_k}\right) +
\sum_{j=0}^{m}   \frac{o(1) }{X_k^j } \right) 
\\
&\geq& 
\sum_{k \in K} \frac{A_k }{ (x_k-x)^{2}} (d_m - o(1)) \max \{ 1, X_k^{-m} \}  > 0. 
\end{eqnarray*}
This proves the first assertion of Proposition \ref{prop1} when $m$ is even, and the case of odd $m$
follows from  Rolle's theorem and the above reasoning applied to $U^{(m+1)}$. 

To prove the second assertion, observe that 
(\ref{ww3}) and (\ref{ww4}) deliver 
$$
U^{(m)}(z) \sim S(z) 
\psi^{(m)}(z) \sim \frac{S(z) m! A_k}{(x_k-z)^{m+1}}
$$
as $z \to x_k$, in which $A_k > 0$ and $S(x)$ has no zeros in $[x_k, x_{k+1}]$. Thus $n_{k,m}$ has the opposite
parity to $m$, and the result follows.
\hfill$\Box$
\vspace{.1in}

Assume henceforth that  
$U^{(m)}$ has finitely many 
non-real zeros. Since all but finitely many zeros and poles of $U$ are real,
\cite[Lemma 2.1]{Lajda18}
and  Proposition \ref{prop1} imply that, as $r \to + \infty$,
\begin{equation}
\label{tsujiest00}
\mathfrak{T}(r, U'/U) = O( \log r ) \quad \hbox{and} \quad 
N(r, 1/U^{(m)}) \leq  2 N(r, U) + O( \log r) \leq (2+o(1)) N(r, U) .
\end{equation} 

\begin{lem}
\label{lemfiniteorder}
$U$ has finite order. 
\end{lem}
\textit{Proof.} Suppose first that $m$ is even.  Proposition \ref{prop1}
gives $k_0 \in \N$ such that if $k \in K$ and 
$|k| \geq k_0$  then $U^{(m)}$ has a simple zero $t_k$ with $x_k < t_k < x_{k+1}$ 
and $x_k t_k > 0$. Moreover,
by assumption and Proposition \ref{prop1}, all but finitely many zeros of
$U^{(m)}$ belong to the set $\{ t_k \}$. Hence the product
$$
\Pi_1(z) = 
\prod_{k \in K, |k| \geq k_0} 
\frac{1-z/t_k}{1-z/x_k} 
$$
converges and maps $H$ into $H$ (by the same argument as in
 Lemma \ref{lemlevost}), and 
\begin{equation}
\label{finiteorderrep}
\frac{U^{(m)}}U = \Pi_1 U_1,
\end{equation} 
where $U_1$ has finitely many zeros and all but finitely many poles of $U$ are
poles of $U_1$. The first estimate of
(\ref{tsujiest00}) and standard properties of the Tsuji characteristic together
lead to $\mathfrak{T}(r, U^{(m)}/U) = O( \log r )$ as $r \to + \infty $. 
Applying  (\ref{C1}) to $\Pi_1$, combined with  the same inequality of Levin and Ostrovskii \cite[p.332]{LeO} as used in Lemma \ref{phitrans}, 
then gives
\begin{eqnarray*}
\frac{T(R, 1/U_1 )}{2 R^2} &\leq& \int_R^\infty \frac{ T(t, 1/U_1)}{t^3} \, dt \\
&\leq& \int_R^\infty \frac{ m(t, 1/U_1)}{t^3} \, dt + O \left( \frac{\log R}{R^2} \right)\\
&\leq& \int_R^\infty \frac{m(t, U/U^{(m)}) + m(t, \Pi_1)}{t^3} \, dt + O \left( \frac{\log R}{R^2} \right)\\
&\leq& 2 \int_R^\infty \frac{\mathfrak{T}(t, U/U^{(m)})}{t^2} \, dt + O \left( \frac{\log R}{R^2} \right)
=  O \left( \frac{\log R}{R} \right)
\end{eqnarray*}
as $R \to + \infty$. Thus $U_1$ has order at most $1$ in the plane, and so,  by (\ref{C1}) applied to $\psi$, 
$$
T(r, U) = m(r, U) + N(r, U) \leq m(r, e^P) + O( \log r ) + N(r, U_1) ,
$$
which implies that the order of $U$ is at most the degree of $P$. 

When $m$ is odd the argument is slightly more complicated. 
In this case, there exists 
$k_0 \in \N$ such that if $k \in K$ and 
$|k| \geq k_0$  then $U^{(m)}$ has in $(x_k, x_{k+1})$ either (a) no zeros at all,
or (b) two 
zeros $u_k, v_k$, these possibly coinciding but having the same sign as $x_k$. 
This time let 
\begin{equation*}
\Pi_1(z) = 
\left( 
\prod
\frac{1-z/u_k }{1-z/x_k} \right) 
\left( 
\prod
\frac{1-z/v_k }{1-z/x_k} \right) 
\end{equation*}
with the products over those $k \in K$ with $|k| \geq k_0$ such that case (b) arises, and
each mapping $H$ into $H$. 
Applying (\ref{C1}) twice then gives 
$m(r, \Pi_1)  = O( \log r) $ as $r \to + \infty$. 
Now define $U_1$ by (\ref{finiteorderrep}): again $U_1$ has finitely many zeros and,
since $m \geq 3$,  all but finitely many poles of $U$ are
poles of $U_1$. 
The remainder of the proof then proceeds as before. 
\hfill$\Box$
\vspace{.1in}

\begin{lem}
 \label{sectorlem}
Let  $K_\varepsilon
 = \{ z \in \C : |z| \geq 1,  \varepsilon \leq | \arg z | \leq \pi - \varepsilon \}$, where $\varepsilon$ is small and positive, and let $n \in \N$. 
Then $U$ satisfies, on $K_\varepsilon$, 
\begin{equation}
 \label{Lasymp}
T_n(z) = \frac{ U^{(n)}(z)}{U(z)} = P'(z)^n (1 + o(1)) \quad \hbox{as $z \to \infty$.}
\end{equation}
\end{lem}
\textit{Proof.}
This is standard and is proved by induction on $n$. 
For $n=1$, (\ref{Lasymp})  is an immediate consequence of (\ref{ww1}) and (\ref{C1}). 
Next,  it may be assumed that $n \geq 1$ and 
 (\ref{Lasymp}) holds on $K_{\varepsilon/2}$, so that (\ref{Lasymp}) for $n+1$ 
 follows from 
Cauchy's estimate for derivatives and the relation
$T_{n+1} = T_n' + T_n T_1$. 
\hfill$\Box$
\vspace{.1in}

\begin{lem}
 \label{lemLmzeros}
Let $\delta, \sigma \in (0, 1)$.
Then $U$  satisfies
\begin{equation}
 \label{polesest}
(m+1 - \delta ) N(r, U) \leq  N(r, 1/U^{(m)})
\end{equation}
as $r \to \infty$ in a set of  lower logarithmic density at least $1 - \sigma$.  
\end{lem}
\textit{Proof.} 
Since $T_m$ has finite order of growth $\rho(T_m)$,
Lemma \ref{Haylem} gives a positive constant $C_1$, depending only
on $\sigma$ and $\rho(T_m)$,  such that
\begin{equation}
 \label{pf4}
 T(2r, 1/T_m) \leq 
T(2r, T_m) + O(1) \leq C_1 T(r, T_m) 
\end{equation}
for  all $r$ in a set $F_1 \subseteq [1, \infty)$ having lower logarithmic density at least 
$1-\sigma$. Let $N_0(r,1/U^{(m)} )$ count common zeros of $U^{(m)}$ and $U$, each
such zero counted only once. Because $S = Re^P$ and $\psi(H) \subseteq H$,   
all but finitely many poles and zeros of $U$ are real, simple and interlaced, and so
\begin{eqnarray*}
(m+1) N(r, U) &=& m N(r, U) + N(r, 1/U) + O( \log r) \\
&\leq& N(r, T_m) + N_0(r,1/U^{(m)}) + O( \log r) \\
&\leq& T(r, 1/T_m) + N_0(r,1/U^{(m)}) + O( \log r) \\
&=& m(r, U/U^{(m)}) + N(r, U/U^{(m)}) + N_0(r,1/U^{(m)}) + O( \log r) \\
&\leq& N(r, 1/U^{(m)})  + m(r, U/U^{(m)}) + O( \log r) 
\end{eqnarray*}
as $r \to + \infty$.  
Now  let 
$\varepsilon$ be small and positive:  then (\ref{Lasymp}) implies that the contribution to 
$m(r, U/U^{(m)})$ from $K_\varepsilon$ 
is bounded as $r \to + \infty$. 
Apply Lemma \ref{lemEF} to $1/T_m = U/U^{(m)} $,
with $R = 2r$ and $\mu(r) = 4 \varepsilon$. In view of (\ref{pf4}) and Lemma \ref{lemfiniteorder}, this
shows that, as $r \to \infty$ in $F_1$, 
\begin{eqnarray*}
(m+1) N(r, U) &\leq& N(r, 1/U^{(m)})  + O( \log r) + 88 \varepsilon \left( 1 +
 \log \frac1{4 \varepsilon } \right)  T(2r, 1/T_m ) \\
  &\leq& N(r, 1/U^{(m)})  + O( \log r) + 88 \varepsilon \left( 1 +
 \log \frac1{4 \varepsilon } \right)  C_1 T(r, T_m ) \\
 &\leq& N(r, 1/U^{(m)})  + O( \log r) + 88 \varepsilon \left( 1 +
 \log \frac1{4 \varepsilon } \right)  C_1 N(r, T_m )  \\
 &\leq& N(r, 1/U^{(m)})  + O( \log r) + 88 \varepsilon \left( 1 +
 \log \frac1{4 \varepsilon } \right)  C_1 (m+1) N(r, U )  ,
 \end{eqnarray*}
again since all but finitely many zeros and poles of $U$ are real, simple and interlaced.
 Because $\varepsilon$ may be chosen arbitrarily small, while $C_1$ 
does not depend on $\varepsilon$, (\ref{polesest}) follows. 
\hfill$\Box$
\vspace{.1in}

To complete the proof of Theorem \ref{thm11}, it remains only to observe that  (\ref{polesest}) 
contradicts (\ref{tsujiest00}),  since $m \geq 2$ and $U$ has infinitely many poles.
\hfill$\Box$
\vspace{.1in}

\section{Proof of Theorem \ref{cor2}} \label{pfcor2}
Let $f$ be as in the hypotheses and assume that $f''$ has finitely many non-real zeros.
Denote by $X$ the set of poles and zeros of
$f$. If $X$ is neither bounded above nor bounded below, using a translation makes it possible to assume that the poles $x_k$ and zeros $y_k$ satisfy
$x_k < y_k < x_{k+1}$ and $x_k/y_k > 0$ for each $k$. Hence $f = \psi e^h$ where $\psi$ is defined as in (\ref{wwpsi}) and maps $H$ into itself, while $h$ is an entire function. If $h$ is constant then it may be assumed that $f(H) \subseteq H$,
so that  (\ref{jdaform}) follows from \cite[Theorem 1.4]{Lajda18} and the remarks preceding it, 
as does  (\ref{hswform})  if $f''$ has only real zeros.
Furthermore, if $h$ is non-constant then a contradiction arises via part (ii) or (iii) of Theorem \ref{thm11}. 

It remains only to consider the case where $X$ is bounded above or below, and here it may be assumed that all zeros and poles of $f$ are positive. If $\min X$ is a pole of $f$ then the argument of the previous paragraph goes through unchanged, and delivers
(\ref{hswform}) or  (\ref{jdaform}), neither of which is compatible with $X$ being bounded below.  
Finally,   if $\min X$ is a zero of $f$ then 
$-1/f = \Psi e^{-h}$, with $h$ entire and $\Psi(H) \subseteq H$: this leads to 
$f = \psi  e^{h}$, where $ \psi = -1/\Psi$ maps $H$ into $H$, 
and the same argument may be deployed. 
\hfill$\Box$
\vspace{.1in}

\textit{Acknowledgement.} The author thanks the referee for a very careful reading of the manuscript and several helpful suggestions.

\noindent
J.K. Langley, Emeritus Professor,\\
Mathematical Sciences, University of Nottingham, NG7 2RD, UK\\
james.langley@nottingham.ac.uk
\end{document}